\documentclass{amsart}
 \usepackage{amsmath}
 \usepackage{epsfig} 
 \usepackage{graphics} 
\usepackage{amsthm}
\usepackage{graphicx,psfrag}

\def\sc{}

\newtheorem*{Theorem0}{\sc Theorem} 
\newtheorem*{MProposition}{\sc Main Proposition}

\newtheorem{Theorem}{\sc Theorem}[section]
\newtheorem{Lemma}[Theorem]{\sc Lemma}
\newtheorem{Proposition}[Theorem]{\sc Proposition}

\newtheorem{Definition}[Theorem]{\sc Definition}

\def\proof{{\bf {\medskip}{\noindent}Proof: }}
\def\proofof#1{{\bf {\medskip}{\noindent}Proof of #1:}}
\def\qed{\hfill$\square$\smallskip}
\def\enddemo{}

\def\intset{\mathop{\hbox{\rm int}}}

\def\det{\mathop{\hbox{\rm det}}}

\def\id{\mathop{\hbox{\rm id}}}

\def\a{\alpha}
\def\b{\beta}
\def\c{\gamma} \def\C{\Gamma}
\def\d{\delta}   \def\D{\Delta}
\def\e{\varepsilon}

\def\l{\lambda}  \def\L{\Lambda}
   \def\O{\Omega}
\def\s{\sigma}
\def\S{\Sigma}
\def\t{\theta}
\def\ttau{\tau}
\def\PPhi{g}

\def\aa{a}
\def\TT{T}
\def\kk{\eta}

\def\Cal{\cal}
\def\cal{\mathcal}
\def\Bbb{\mathbb}
\def\dd{\delta}

\def\disp{\displaystyle}

\title [A hyperbolic diffeomorphism near identity]
{A hyperbolic diffeomorphism with countably many ergodic
components near identity}
\author[Huyi Hu and Anna Talitskaya]{}


 \email{hu@math.msu.edu}
 \email{anjuta@\@math.northwestern.edu}

\begin{document}
\maketitle

\centerline{\scshape   Huyi Hu and Anna Talitskaya}
\medskip

  {\footnotesize 
  \centerline{ Huyi Hu: Department of Mathematics,}
  \centerline{ Michigan State University}
  \centerline{ East Lansing, MI 48824}
  \centerline{ Anna Talitskaya: Department of Mathematics }
  \centerline{ Northwestern  University} 
  \centerline{ Evanston, IL 60201, USA} }
\medskip

\medskip

\begin{abstract}
We construct a smooth hyperbolic volume preserving diffeomorphism  
on a four dimensional compact Riemannian manifold  
which has countably many ergodic components 
and  is arbitrarily close to the identity map.
\end{abstract}


\section{Introduction}
\setcounter{equation}{0}
In 1978 Ya. Pesin proved that a hyperbolic volume preserving
diffeomorphism $f$ on a compact Riemannian manifold $M$ admits 
at most countably many ergodic components (mod 0) by showing that 
almost every ergodic component has positive volume (\cite{P}).
Here hyperbolicity means nonzero Lyapunov exponents almost everywhere.
However, an example of a system which has infinitely many
ergodic components was unknown till 2000 
when D. Doglopyat, Ya. Pesin and H. Hu constructed a diffeomorphism 
of the 3-dimensional torus that satisfies the properties stated 
above (see \cite{DHP}).

The construction starts with a  product $F$ of a two dimensional 
toral automorphism $A:\mathbf{T^2}\to \mathbf{T^2}$ and the identity 
map on a circle, $F=A\times id$. 
The space is then partitioned into countably many subsets, 
$M=\mathbf{T^2}\times S^1=\cup_{i=0}^\infty  
\big(\mathbf{T^2}\times [ \frac{1}{2^{i+1}},\frac{1}{2^i} ]\big)$
The diffeomorphism  is then obtained
by applying two small perturbations on each subset,  one for
ergodicity, and the other one for hyperbolicity.  The latter one is
a perturbation which removes the zero Lyapunov exponent along 
the interval direction.
Since the original map $F$ is a product of Anosov and
identity, the resulting map will be away from the identity. 
In this paper we construct a system satisfying the properties in \cite{DHP} 
in the proximity of the identity map. Our map is a perturbed 
 geodesic flow on a surface of constant negative curvature
instead of a toral automorphism.  We will need three perturbations 
instead of two since there are two zero Lyapunov exponents to remove.

There are two difficulties in this construction.  First, to obtain 
ergodicity we need accessibility.  For this purpose we need to show
that almost every point is accessible to a given central leaf
along the interval direction.  Since both strong stable and unstable 
leaves are one dimensional, the accessibility is easy to achieve 
in a three dimensional space as in \cite{DHP} and  extra work is 
necessary in our case because the space we work with is four
dimensional.  Second, the technique to remove 
the second zero Lyapunov exponent is more delicate.
While we change the last Lyapunov exponent, we need to keep all other 
exponents nonzero at the same time.
But the direction corresponding to the ``small'' nonzero exponent obtained 
by the earlier perturbation is not stable under perturbation. 

Our results show that systems which are hyperbolic, volume preserving,
and with countably many ergodic components exist in any neighborhood
of the identity.
However, the phenomenon may not be true for other properties
of dynamical systems. 
For example, it is generally believed that there is no uniformly
hyperbolic (Anosov) system in a neighborhood of the identity.

It is also interesting to compare our theorem with the results 
in \cite{CS} and \cite{X}.  They show that on any manifold $M$
of dimensional at least two, there are open sets of volume
preserving diffeomorphisms of $M$,  all of which have positive 
measure sets on which all of the Lyapunov exponents are zero.
In the terminology used in \cite{BPSW}, they prove that there 
are open sets of {\it neutral} volume preserving diffeomorphisms
{\it away} from the identity, 
while we show that there are {\it (fully) nonuniformly hyperbolic} 
systems {\it near} identity.

\section{Statement of Results}

We prove the following result.

\begin{Theorem0}
There exists a four dimensional compact Riemannian manifold $M$
such that for any $\dd_0>0$, we can find a $C^\infty$ diffeomorphism
$f$ of $M$ with the following properties:
\begin{enumerate}
\item $\|f-\id\|_{C_1}\le \dd_0$;
\item $f$ preserves the Riemannian volume $\mu$ on $M$;
\item $\mu$ is a hyperbolic measure;
\item $f$ has countably many ergodic components which are open $\pmod 0$.
\end{enumerate}
\end{Theorem0}

\section{Construction}

Let $\PPhi^t: M_0\to M_0$ be a geodesic flow on a compact surface
of a negative constant curvature.

Choose a closed orbit ${\cal C}$.

Take $0 < \delta\le \delta_0/2$ such that all points in ${\cal C}$
are periodic points of $\PPhi^\delta$.  Take $G=\PPhi^\delta$.

Let $F=G\times \id$ be the map from $M=M_0\times \Bbb S^1$ to itself.
We will perturb $F$ to obtain the desired map $f$.

Consider a countable collection of intervals $\{I_n\}_{n=1}^{\infty}$
on the circle $\Bbb S^1$, where
$$
I_{2n}=[(n+2)^{-1}, (n+1)^{-1}], \quad
I_{2n-1}=[1-(n+1)^{-1},1-(n+2)^{-1}].
$$
Clearly, $\bigcup\limits_{n=1}^\infty I_n=(0,1)$ and
$\intset I_n$ are pairwise disjoint.

Choose $\d'>0$ small enough.
By Main Proposition below, for each $n$ one can construct
a $C^\infty$ volume preserving ergodic hyperbolic diffeomorphism
$f_n: M_0 \times I\to M_0\times I$
satisfying: 1) $\|F-f_n\|_{C_1}\le \d' n^{-4}$; \      
2) for all $0\leq n<\infty$,
$D^k f_n|_{M_0 \times\{s\}}=D^k F|_{M_0 \times\{s\}}$
for $s=0$ or $1$.

Let $L_n: I_n\to I$ be the affine map and
$\pi_n=(\id, L_n): M_0\times I_n\to M_0\times I$.
Clearly, $\|\pi_n\|\le 5n^2$ and $\|\pi_n^{-1}\|\le 1$.

We define the map $f$ by setting
$f|_{M_0\times I_n}=\pi_n^{-1}f_n\pi_n$ for all $n$
and then letting $f|_{M_0\times \{0\}}=F|_{M_0\times \{0\}}$.
Note that
$$
\bigl\|F|_{M_0\times I_n}-\pi_n^{-1}f_n\pi_n\bigr\|_{C_1}
\le \bigl\|\pi_n^{-1}\bigl(F -f_n\bigr)\pi_n\bigr\|_{C_1}
\le \d' n^{-4}\cdot 5n^2
= 5\d' n^{-2}\le \dd.$$
It follows that $f$ is $C^\infty$ on $M$ and has the required
properties.

\section{Main Proposition}
\setcounter{equation}{0}

The goal of this section is to prove the following statement.
\begin{MProposition}
Let $S=G\times \id$ be the diffeomorphism
from $N=M_0\times I$ to itself.
For any  $\d>0$, there exists a map $P$ such that:
\begin{enumerate}
\item  $P$ is a $C^\infty$ volume preserving diffeomorphism of $N$;
\item $\|S-P\|_{C_1}\le \d$;
\item  for all $0\leq n<\infty$,
$D^nP|_{M_0\times\{s\}}=D^nS|_{M_0\times\{s\}}$ if $s=0$ and $1$;
\item $P$ is ergodic with respect to the Riemannian volume and has
non-zero Lyapunov exponents almost everywhere.
\end{enumerate}
\end{MProposition}

Note that $S$ is not ergodic, and has two zero Lyapunov exponents.
We perturb $S$ by three small perturbations 
$h^{(i)}: \Omega_i\to \Omega_i$, $i=1,2,3$, where $\Omega_i\subset N$,
to get ergodicity and to remove zero Lyapunov exponents.

\proofof{Main Proposition}
Note that the tangent bundle of $N$ can be written
as a direct sum of four one-dimensional $S$-invariant subbundles:
$$
TN=E^u(S)\oplus E^s(S) \oplus E^c(S) \oplus E^n(S),
$$
where $E^u(S)$, $E^s(S)$ and $E^c(S)$ are the unstable, stable and
flow directions of the geodesic flow $\PPhi^t: M_0\to M_0$, and
$E^n(S)$ is the tangent space of $I$. The corresponding Lyapunov
exponents of $S$ at $w\in N$ are denoted by $\l^u(w, S)$, $\l^s(w,
S)$, $\l^c(w, S)$ and $\l^n(w, S)$ respectively. It is clear that
the Lyapunov exponents are constants at almost every point $w\in
N$, though $S$ is not ergodic. We simply denote them by $\l^u(S)$
etc. Also, we know that $\l^u(S)>0$, $\l^s(S)<0$ and
$\l^c(S)=\l^n(S)=0$. We will often take local coordinate system
$w=(x,y,t,z)$ in $N$ in such a way that
\begin{eqnarray}\label{f3.0}
d\mu=dxdydtdz, \ \mbox{where}\
{\partial \over \partial x}=E^u(S), \
{\partial \over \partial t}=E^c(S), \
{\partial \over \partial z}=E^n(S).
\end{eqnarray}


Fix $\tau\in (0, 2/3)$.

Take $ k_0\in \Bbb N$ such that for $\c\le \d$, $\tau_2=1/4$, 
some $C>0$, and $\theta=\pi/ k_0$ such that the properties stated in
Lemma~\ref{L6.4} are satisfied.

Recall that ${\cal C}$ is the closed orbit taken in Section 2.
Take another closed orbit ${\cal C}'\subset M_0$ of $g^t$.

Choose a set $\O_0$ and constant $\e_0$ according to Lemma~\ref{L7.1} 
with $\tau_1=(1-1.3\tau)(k_0+1)^{-1}$.   Hence
\begin{eqnarray}\label{f3.3}
\mu (\bigcup_{0\le i\le k_0}G^i \O_0)\le 1-1.3\tau.
\end{eqnarray}

We also assume that $\e_0$ is small enough such that 
${\cal C}$ and ${\cal C}'$ are at least $3\e_0$ separated.

\smallskip
{\it Construction of $h^{(1)}: \Omega_1\to \Omega_1$}.

Fix $p\in {\cal C}$, a periodic point of $G$ with period $m\in \Bbb N$.

We assume further that $\e_0>0$ is small such that for any $w\in N$,
\begin{eqnarray}\label{f3.1}
\mu (B(w, \e_0 ))< \frac{0.1\tau}{k_0 m}
\end{eqnarray}
Here $B(w,\e_0)$ denotes the ball in $M_0$ of radius $\e_0$
centered at $w$.

We assume that $G^i(B(p,\e_0 ))\cap B(p,\e_0 )=\emptyset$ for
$i=1,\cdots m-1$.

For any $z\in M_0$, let $V^u(z)$ and $V^s(z)$ be the local unstable and
stable one dimensional manifold at $z$ for $G$ of ``size'' $\e_0$.

Since both $W^{u}({\cal C}')$ and $W^{s}({\cal C}')$
are dense in $M_0$, we can choose $p_1, p_2\in {\cal C}'$,
and the smallest integers $n_1, n_2>0$ such that each intersection
$$
G^{-n_1}V^s(G^{n_1}p_1)\cap V^u(p)\cap B(p,\e_0) \quad
\mbox{and}\quad
G^{n_2}V^u(G^{-n_2}p_2)\cap V^s(p)\cap B(p,\e_0)
$$
consists of a single point $q_1$ and $q_2$ respectively.

Take $\e_1 \le \min\{\d, d(p,q_1)/2, d(p,q_2)/2\}$.
Take $\ell\ge 2$ such that
\begin{eqnarray}\label{f3.2}
G^{-\ell m} (q_1)\not\in B(p,\e_1),\quad
G^{-(\ell+1)m}(q_1)\in B(p,\e_1). \qquad
\end{eqnarray}
Then we take $\e_2\in (0,\e_1)$ such that $G^{-(\ell+1)m}(q_1)\in B(p,\e_2)$.

\begin{figure}
\begin{center}
\psfrag{C}{$ C$}
\psfrag{C1}{$C'$}

\psfrag{B(p,e0)}{$B(p,\epsilon_0 )$}

\psfrag{p1}{$p_1$}

\psfrag{p2}{$p_2$}

\psfrag{B(p,e1)}{$B(p,\epsilon_1 )$}
\psfrag{q2}{$q_2$}
\psfrag{q1}{$q_1$}

\psfrag{B(p,e2)}{$B(p,\epsilon_2 )$}

\psfrag{Vu(p1)}{$V^u (p_1 )$}

\psfrag{Vs(p2)}{$V^s (p_2 )$}

\psfrag{GVs(Gp1)}{$G^{-n_1}V^s (G^{n_1}p_1 )$}

\psfrag{GVu(p2)}{$G^{n_2}V^s (G^{-n_2}p_2 )$}

\psfrag{p}{$p$}

\psfrag{Vs(p)}{$V^s (p)$}
\psfrag{Vu(p)}{$V^u (p)$}

\includegraphics{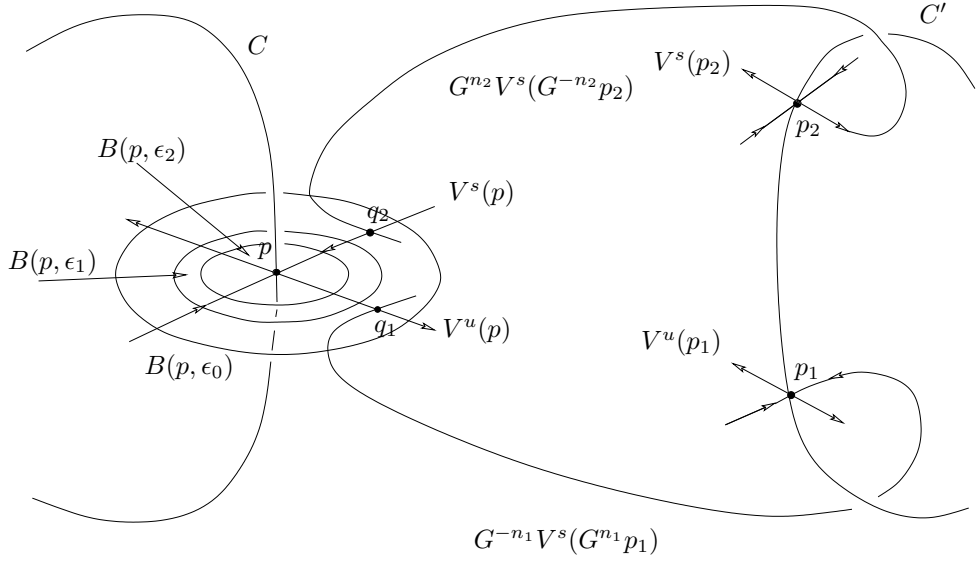}
\end{center}
\caption{Construction of $h^{(1)}$.}
\end{figure}

Let $\Omega_1=B(p,\e_0 )\times I$.  Denote 
\begin{eqnarray}\label{f3.3a}
\tilde \O_1=\O_1\bigcup \bigl( B({\cal C} ,\e_3 )\times I\bigl) \bigcup
\bigl( \cup^{m-1}_{i=1} G^i (B(p,\e_0 ))\times I\bigl) \bigcup 
\O '_1 \bigcup \O ''_1,
\end{eqnarray}
where 
\begin{eqnarray*}
&&\O '_1 
= B\bigl(\cup^{\infty}_{i=0} G^{-n_1 +i}V^s (G^{n_1}p_1),\e_3\bigr)
  \times I,                  \nonumber \\
&&\O ''_1 
  =B\bigl(\cup^{\infty}_{i=0} G^{n_2 -i}V^u (G^{-n_2}p_2),\e_3\bigr)
\times I,
\end{eqnarray*}
and $B(\O ,\e )$ is the $\e$-neighborhood of the set $\O$ in $M_0$ and
$\e_3$ is chosen such that 
\begin{eqnarray}\label{f3.4}
\mu \tilde \O_1\le 0.1\tau /k_0.  
\end{eqnarray}
This is possible because of $(\ref{f3.1})$.



Take a coordinate system $w=(x,y,t,z)$ in
$\Omega_1$ satisfying (\ref{f3.0}). Take the $y$-coordinate 
in such a way that
$\partial\over\partial y$ $ = E^s (S)$ along the path $V^s (p)$.

Choose a $C^\infty$ function $\phi=\phi(r): \Bbb R^+\to \Bbb R^+$
which satisfies
\begin{enumerate}
\item $\phi(r)=\phi_0$ if $r\in [0, \e_2 ]$,
where $\phi_0$ is a positive constants;
\item $\phi(x)=0$ if $x\ge \e_1$;
\item $\phi '(x)\le 0$ for any $x$;
\item $\|\phi\|_{C^1}\le\d$.
\end{enumerate}

Choose a $C^\infty$ function $\psi=\psi(y): \Bbb R\to \Bbb R^+$
which satisfies
\begin{enumerate}
\item[5.] $\psi(x)=\psi_0$ if $x\in (-\e_2, \e_2)$, where $\psi_0$
is a positive constant;
\item[6.] $\psi(x)=0$ if $|x|\ge \e_1$;
\item[7.] $\|\psi\|_{C^1}\le\d$;
\item[8.] $\int_0^{\pm\e_1}\psi(s)ds=0$.
\end{enumerate}

We also choose a $C^\infty$ functions $\xi: I\to \Bbb R^+$ satisfying:
\begin{enumerate}
\item[9.] $\xi(s)> 0$ on $(0,1)$;
\item[10.] $\xi^{(i)}(0)=\xi^{(i)}(1)=0$ for $i=0,1,2,\cdots$;
\item[11.] $\|\xi\|_{C^1}\le \d$.
\end{enumerate}

Then we define the vector field $X$ on $\O_1$ by
$$
X(x,y,t,z)
=\Bigl(-\phi (\sqrt{y^2+t^2} )\xi'(z)\int_0^x\psi(u)du,
\ \ 0, \ \ 0, \ \
\phi (\sqrt{y^2+t^2} )\xi(z)\psi(x) \ \ \Bigr).
$$
It is easy to check that $X$ is a divergence free vector field
supported in
$(-\e_1, \e_1)\times (-\e_1, \e_1)\times(-\e_1, \e_1)\times I\in\O_1$.
We define the map $h^{(1)}=h^{(1)}_\b$ on $\O_1$ to be the time $\b$ map
of the flow  generated by $X$ and we set $h^{(1)}=\id$ on the
complement of $\O_1$. It is easy to see that $h^{(1)}$ is a $C^\infty$
volume preserving diffeomorphism.
Also, we assume that $\b$ is small enough such that
$\|h^{(1)}-\id\|\le \d$.

Let $R=R_\b=h^{(1)}_\b\circ S$ for some small $\b>0$. By
Proposition~\ref{P4.1}, for any small $\b>0$, $R: N\to N$ is ergodic.

\smallskip
{\it Construction of $h^{(2)}: \Omega_2\to \Omega_2$}.

Choose a point $w^*=(w^*_0, z_0)\in N$, where $w^*_0\in M_0$ and $z_0$
is the midpoint of $I$.  Choose $0<\e_4\le\e_0$.  
Let $\O_2=B(w^*, \e_4)$, where $B(w,\e)$ denotes the ball in $N$ 
of radius $\e$ centered at $w\in N$.
We assume that $w^*_0$ and $\e_4$ is chosen in such a way that
$\O_2\cap (\tilde\O_1 \cup (\O _0 \times I))=\emptyset$ and
$S^{-i}\O_2\cap (\O_2\cup \O_1 )=\emptyset$
for $i=1, \cdots, N_0$, where $N_0$ is an integer large enough such that
(B. 19) in \cite{DHP} holds for this $N_0$.
We also assume that $\e_4$ is small enough such that
\begin{eqnarray}\label{f3.6}
\mu (\O_2 ) \le \frac{0.1\ttau}{ k_0}.
\end{eqnarray}

Let $\rho(r)=(\e_4/\e_1)\phi(r\e _1 /\e_4)$.
To define $h^{(2)}$, we take the cylindrical coordinate system
$(r,\t, y, t)$ in $\Omega_2 $, where $x=r\cos\t$, $y=y$, $t=t$
and $z=r\sin\t$.
Define
$h^{(2)} = h^{(2)}_{\a}$ on $N$ by
\begin{eqnarray}\label{f3.5}
 h^{(2)}(r,\t, y, t)
=(r, \  \t+\a\rho(\sqrt{y^2+t^2})\rho(r),\  y, t).
\end{eqnarray}
Then we extend $h^{(2)}$ to $N$ by letting $h^{(2)}=\id$ on $\O_2^c$.

Let
$Q=Q_{\a\b}=R_\b \circ h^{(2)}_\a=h^{(1)}_\b\circ S\circ h^{(2)}_\a$.

We denote by $\kk_{\alpha\beta}(w)$ the expanding rate of
$Q_{\alpha\beta}$ along its unstable direction
$E^u_{w}(Q_{\alpha\beta})$. By Proposition~\ref{P5.1},  
we have $\int_N \log \kk_{\alpha 0}(w)dw<\lambda^u (S)$. 
Since $\kk_{\alpha\beta}$
change smoothly with $\beta$ we conclude that for some small
$\beta >0$
$$
\int_N \log \kk_{\alpha\beta}(w)dw <\lambda^u (S).
$$

By Proposition~\ref{P4.1}, $Q_{\alpha\beta}$ is ergodic. 
So we denote by $\l^u(Q_{\alpha\beta})$ the largest Lyapunov 
exponent of $Q_{\alpha\beta}$.  By the above inequality,
$\l^u(Q)\le \lambda^u (S)$.

Note that both $Dh^{(1)}$ and $Dh^{(2)}$ preserve $E^{un}(S)$
bundle, and  for any $w\in N$, $|\det Dh^{(1)}_w|_{E^{un}(S)}|=|\det
Dh^{(2)}_w|_{E^{un}(S)}|=1$. 
So $DQ$ preserves $E^{un}(S)$ bundle, and 
$|\det DQ_w|_{E^{un}(S)}|=|\det DS_w|_{E^{un}(S)}|$. 
Let $\l^n(Q)$ denote the other Lyapunov exponent on $E^{un}(S)$. 
Then we have
$$
\l^u(Q)+\l^n(Q)=\l^u(S)+\l^n(S).
$$
Since $\l^u(Q)< \l^u(S)$ and $\l^n(S)=0$, we have $\l^n(Q)>0$.

Note that the perturbations also preserve $E^{ucn}(S)$ bundle.
So we have $$\l^u(Q)+\l^c(Q)+\l^n(Q)=\l^u(S)+\l^c(S)+\l^n(S)$$
and therefore $\l^c(Q)=0$, where $\l^c(Q)$ is the third Lyapunov
exponent of $Q$ on $E^{ucn}(S)$.
Applying the same arguments, we also get $\l^s(Q)=\l^s(S)< 0$, though
the stable bundle $E^s(Q)$ may not be equal to $E^s(S)$.

\smallskip
{\it Construction of $h^{(3)}: \Omega_3\to \Omega_3$}.

Denote $\l=\l^n(Q)$. We assume that $\l$ is small in comparison
with  $\l ^u (Q)$.

Denote
$$
\Lambda'=\Lambda'(K)
=\{w\in N: |\log |DQ_w^k|_{E^n(w, Q)}|-k\lambda|\le 0.1k\lambda,
\forall |k|\ge 0.5K\},
$$
and  $\disp \Lambda=\bigcap_{i=0}^{k_0}Q^{-i}\Lambda'$.
Note that $\mu \Lambda'\to 1$ as $K\to \infty$.
We assume that $K$ satisfies
\begin{eqnarray}\label{f3.7}
&K\l\ge \max\{2k_0\l, \ 1.25\log 2, \ -10k_0\log(1-\d)\},  \\
\label{f3.8}
&0.001\tau^2\l+ \mu \L^c \log(1-\d)>0,  \\
\label{f3.9}
&\mu \L^c \le 0.1\ttau,
\end{eqnarray}
where $\L^c$ is the complement of $\L(K)$ in $N$.

Note that if $w\in\L'$ then
$\|DQ^n_w (v)\| \ge e^{0.9 n\l }\| v\|$ for
$n\ge 0.5K$, $v\in E^{un}_w (Q)$.

Denote
$$
\O=\Lambda^c\bigcup
\Bigl(\bigcup_{i=0}^{k_0}Q^{-i}
\bigl(\O_0\cup \tilde \O_1\cup \O_2\bigr)\Bigr).
$$
By (\ref{f3.9}), (\ref{f3.3}), (\ref{f3.4}) amd (\ref{f3.6}),
\begin{eqnarray}\label{f3.10}
\mu \O^c\ge 1-0.1\ttau -(1-1.3\tau )-0.1\tau -0.1\tau =\tau.
\end{eqnarray}

Choose a set $\C'\in N$ such that $Q^i\C'\cap \C'=\emptyset$ ,
$ -K\le i\le 5\tau^{-1}K-2K-1$, $i\not=0$.
Here we assume that $\tau^{-1}$ is an integer, otherwise
we can use a smaller $\tau$ instead.  Denote
$$
\disp \overline \C'=\bigcup_{i=-K}^{5\tau^{-1}K-2K-1}Q^i\C'.
$$
We also require that $\mu \overline \C'$ is close to $1$ such that
\begin{eqnarray}\label{f3.11}
(1-0.5\tau )\cdot\mu \overline \C'\ge 1-0.6\tau.
\end{eqnarray}
The choice of such $\C'$ is possible
because of the Rokhlin-Halmos Lemma.

We define
$$\C_0=
\{Q^j w: \  w\in\C', 0\le j \le 4\tau ^{-1}K-K-k_0,
            Q^j w\in\O^c , Q^i w\notin\O^c \mbox{ for } i\le j\}.
$$
In other words, $\C_0$ is the set of points from each trajectory
$\{ Q^i w\}_{i=0}^{4\tau^{-1}K-K-k_0}$ that enter the set $\O^c$
the first time.


Clearly $Q^i\C$, $i=-K, \cdots, K+k_0$, are pairwise disjoint.
Let $\C_i=Q^i\C$ for $i=-K, \cdots, K+k_0$,
$\disp \C_{j,k}=\cup_{i=j}^k\C_i$ for $j\le k$,
and in particular, $\overline \C=\C_{-K, K+k_0}$.
Since $\C$ is disjoint with $\O$, it is clear that
$\C_i\cap(\O_0\cup \tilde \O_1\cup \O_2)=\emptyset$ 
for $i=1, \cdots, k_0$.
Since we are going to make perturbation $h^{(3)}$
 around the set $\C_{0,k_0 -1}$,
this condition guarantees that the properties of $Q$ we mentioned above
still remain.

Approximate $\C=\C_0$ by finitely many number of disjoint sets
of the form
$$
\D_{0j}=B^u(x_j, r_j')\times B^s(y_j, r_j'')\times B^{cn}((t_j,z_j), r_j),
$$
where $w_i=(x_j, y_j, t_j, z_j)\in N$, $r_j', r_j''\ge r_j$ for
$j=1, \cdots, J$ and $B^u$, $B^s$ and $B^{cn}$ are the balls  that
correspond to the $x$, $y$ and $zt$-coordinates.
Denote $\D_{ij}=Q^i\D_{0j}$, $\disp\D_i=\cup_{j=1}^J \D_{ij}$ for
$i=-K,\cdots, K +k_0$. We can pick $\D_{0j}$ in such a way that
$ \D_{ij}\cap \D_{kl}=\emptyset$ for $(i,j)
\ne (k,l)$, $-K\le i,k\le K+k_0$, $1\le j,l\le J$
 and $\D_{ij}\cap (\O_0\cup \tilde\O_1\cup \O_2) =\emptyset$ 
for $0\le i\le k_0$, $0\le j\le J$.

Also, denote $\disp\overline \D=\cup_{i=0}^{k_0 -1}\D_i$.
Clearly, $\D_i$ is an approximation of $\C_i$ for $i=1, \cdots, k_0$.
We may assume that for each $i=0, \cdots, k_0$,
\begin{eqnarray}\label{f3.12}
\mu(\C_i\triangle \Delta_i)\le 0.05\max\{\mu\C_i, \mu\D_i\}.
\end{eqnarray}

Take $\O_3=\overline \D$.
On each $\Delta_{ij}$ we apply Lemma~\ref{L6.4} to get a map $h=h_{ij}$
and a subset $\D_{ij}'\subset \D_{ij}$
such that $\|h_{ij}-\id\|\le \d$, $\mu \D_{ij}'/\mu\D_{ij}\ge 3/4$
and restricted to $\D_{ij}'$, $h$ is a rotation of angle $\pi/2k_0$
along the $cn$ plane.   Note that $DS|_{E^{cn}}=\id$,
we can require that $Q\D_{ij}'=\D_{i+1,j}'$
for $j=0, \cdots, k_0-1$.

Let $\D_{i}'=\cup_{j=1}^J \D_{ij}'$.

Let $h^{(3)}=h_{ij}$ on each $\D_{ij}$ and $h^{(3)}=\id$
otherwise. Take $P=Q\circ h^{(3)}$.
By Proposition~\ref{P4.1}, $P$
is ergodic.
By Proposition 6.1, we have that for almost every $w\in N$, for
all $v\in E^{ucn}(w, S)$, $\lambda(w, v, P)>0$. So $P$ has three
positive Lyapunov exponents. Since the foliation $W^{ucn}(S)$ is
preserved, we have $E^{ucn}(w, S)=E^{ucn}(w, P)$.  This implies
$\l^s(P)=\l^s(S)<0$.
\qed \enddemo

\section{Ergodicity}
\setcounter{equation}{0}

Recall that for a diffeomorphism $f:M\to M$, two points $w_1, w_2\in M$
are called {\it accessible (with respect to $f$)}
if they can be joined by a piecewise differentiable piecewise nonsingular
path which consists of segments tangent to either $E^u(f)$ or $E^s(f)$.
The diffeomorphism $f$ is {\it essentially accessible}
if almost any two points in $M$ (with respect to the Riemannian volume)
are accessible.


\begin{Proposition}\label{P4.1}
For any $\a, \b, \c>0$ sufficiently small,
the diffeomorphisms $R=R_\a$, $Q=Q_{\a\b}$ and $P=P_{\a\b\c}$
of $N$ are ergodic with respect to the Riemannian volume.
\end{Proposition}

\proof
By a result of Pugh and Shub, if a $C^2$ diffeomorphism is
partially hyperbolic, center bunched, dynamically coherent
and essentially accessible, then the diffeomorphism is
ergodic.  (See \cite{PS}, Theorem A, also see \cite{BPSW}, Theorem 2.2.)

Clearly, all $S, R, Q$ and $P$ are partially hyperbolic and
center bunched if $\a, \b$ and $\c$ are small.

Note that $S$ is dynamically coherent and its center foliation
is plague expansive.   Since $R, Q$ and $P$ are $C^1$ close
to $S$, they are all dynamically coherent if $\a, \b$ and $\c$ are
small by a theorem of Pugh and Shub (\cite{PS}, Theorem 2.3).

By the lemma below, $R, Q$ and $P$ are essentially accessible.
So they are ergodic.
\qed 
\enddemo


\begin{Lemma}\label{L4.2}
Any two points $w, w'\in \intset N$ are accessible with respect
$R, Q$ or $P$.  Therefore,  $R, Q$ or $P$ are essentially accessible.
\end{Lemma}

\proof
Recall that $p\in {\Cal C}$ is chosen to construct $h^{(1)}$.
Denote $I_p=\{p\}\times (0,1)$, and $\bar I_p=\{p\}\times [0,1]$.
For $T=R, Q$ or $P$, let ${\Cal A}_T$ be the set of points
that are accessible to some point in $\bar I_p$ with respect
to $T$.  Since ${\Cal A}_T$ is both open and closed
in $N$, we get ${\Cal A}_T=N$.  Further, since $\partial N$
is unperturbed by any of $h^{(i)}$, $i=1,2,3$,
any point in $\intset N$ is accessible to a point in $I_p$.

By the Sublemma below, we know that any two points in $I_p$
are accessible.  Since the accessibility property is
symmetric and transitive, we get that any two points
in $\intset N$ are accessible.
\qed \enddemo

\begin{Lemma}\label{L4.3}
Let $T=R, Q$ or $P$. For any $s\in (0,1)$,
\begin{eqnarray}\label{f4.1}
{\Cal A}_T(p,s)\supset I_p,
\end{eqnarray}
where ${\Cal A}_T(p,s)$ is the set of points accessible to
$(p,s)\in N$ with respect to $T$.
\end{Lemma}

\proof
In the arguments below we only use the sets $\tilde \O_1$
defined in (\ref{f3.3a}), $\O_0$ chosen in Lemma~\ref{L7.1},
and the fact that the strong stable and unstable foliations
are continuous with the maps.  Since $\O_2$ and $\O_3$ are disjoint
with $\O_0$ and $\tilde \O_1$, and $h^{(i)}$ are identity
outside $\O_i$ for $i=2, 3$, the arguments work for all
$R, Q$ and $P$.  So we drop the subscript $T$ in ${\Cal A}_T(p,s)$
and simply write ${\Cal A}(p,s)$ instead.

We use the coordinate system $(x,y,t,z)$ in $\Omega_1$ described
as we construct $h^{(1)}$.
Denote $h=h^{(1)}_\a$.
Since the map $h$ preserves the leaf $I_p$, we have that
$$
h(\bar 0,z)=(h^1(\bar 0,z),\,h^2(\bar 0,z),\,h^3(\bar 0,z),\,
h^4(\bar 0,z)) =(\bar 0,h^4(\bar 0,z))
$$
for all $z\in (0,1)$, where $\bar 0=(0,0,0)$.
It suffices to show that for every $z\in (0,1)$,
\begin{eqnarray}\label{f4.2}
{\Cal A}(p,z)\supset \{(p,z'): z'\in [(h^{-\ell}_t)^4(p,z), z]\},
\end{eqnarray}
where $\ell$ is chosen by (\ref{f3.2}). In fact, since accessibility is
a transitive relation and $h^{-n}_t(p,z)\to (p,0)$ for any $z\in (0,1)$
as $n\to \infty$,
(\ref{f4.2}) implies that ${\Cal A}(p,z)\supset \{(p,z'): z'\in (0,z]\}$.
Since this holds true for all $z\in (0,1)$
and accessibility is a reflective relation, we obtain (\ref{f4.1}).

Now we proceed with the proof of (\ref{f4.2}).

Recall that $p_1,p_2\in {\cal C}'$,
$q_1\in G^{-n_1}V^s(p_1)\cap V^u(p)$ and
$q_2\in G^{n_2}V^u(p_2)\cap V^s(p)$.
For $z_0=z$, we choose $z_i$, $i=1,\cdots,5$, such that
\begin{eqnarray*}
&(q_1, z_1)\in V^u((p,z_0), \TT),   \qquad
G^{-n_1}V^s((p_1, z_2), \TT)\ni (q_1, z_1),  \\
&(p_2, z_3) \  \mbox{is accessible to} \  (p_1, z_2), \\
&(q_2, z_4)\in G^{n_2}V^u((p_2, z_3), \TT), \qquad
V^s((p, z_5), \TT)\ni (q_2, z_4).
\end{eqnarray*}
This means that $(p,z_5)\in{\Cal A}(p,z_0)$.

Let $\pi: N=M_0\times I\to M_0$ be a projection.
Note that $\pi V^s(w, \TT)=V^s(\pi w, G)$ if $w\in \{p\}\times I$ 
or $G^{-n_1}V^s((p_1, G)\times I$, and
$\pi V^u(w, \TT)=V^u(\pi w, G)$ if $w\in \{p\}\times I$
or $G^{n_2}V^u((p_2, G)\times I$.
In other words, $V^s(\pi w, \TT)\in V^s(w, G)\times I$
and $V^u(\pi w, \TT)\in V^u(w, G)\times I$.
Hence we know that $z_1, z_2, z_4, z_5$ are uniquely determined by
$z_0, z_1, z_3, z_4$ respectively.
We will show that
\begin{enumerate}
\item  $z_1=(h^{-\ell})^4(p,z_0)$;
\item  $z_2=z_1$;
\item  $z_3$ can be chosen arbitrarily close to $z_2$;
\item  $z_4=z_3$;
\item  $z_5\le z_4$.
\end{enumerate}
If so, we can get
\begin{eqnarray}\label{f4.3}
z_5\le (h_t^{-\ell})^4(p,z_0).
\end{eqnarray}
By continuity, we conclude that
$$
\{(p,z'): z'\in [z_5, z_0]\}\subset {\Cal A}(p,z_0)
$$
and therefore (\ref{f4.2}) follows.

By the construction of $h=h^{(1)}_\a$, we know that for any
$q\in B(p,\e_2)$,
\begin{eqnarray}\label{f4.4}
h^4(q,z)=h^4(p,z).
\end{eqnarray}
Recall that $p$ is a periodic point of $G$ with period $m$.
We have
$$
\TT^{-m}(p,z_0)=h^{-1}S^{-m}(p,z_0)=h^{-1}(G^{-m}p, z_0)
=(p,(h^{-1})^4(p,z_0) ),
$$
and therefore for any $k\ge 1$,
\begin{eqnarray}\label{f4.5}
\TT^{-km}(p,z_0)=(p, (h^{-k})^4(p,z_0) ).
\end{eqnarray}
On the other hand, $G^{-km}q_1\notin B(p,\e)$ if $k\le \ell$.
It follows that
$\TT^{-km}(q_1,z_1)=S^{-km}(q_1,z_1)=(G^{-km}q_1, z_1)$.
Hence, by (\ref{f3.2}),
\begin{eqnarray*}
&\TT^{-(\ell+1)m}(q_1,z_1)=h^{-1}S^{-(\ell+1)m}(q_1,z_1) \\
=&h^{-1}(G^{-(\ell+1)m}q_1, z_1)=(q^{(1)}, (h^{-1})^4(p,z_1) )
\end{eqnarray*}
for some $q^{(1)}\in B(p,\e_2)$ and therefore, for any $k>0$,
\begin{eqnarray}\label{f4.6}
\TT^{-(\ell+k)m}(q_1,z_1)=(q^{(k)}, (h^{-k})^4(p,z_1) )
\end{eqnarray}
for some $q^{(k)}\in B(p,\e_2)$.  
Since $(q_1, z_1)\in V^u((p,z_0), \TT)$, 
by (\ref{f4.5}) and (\ref{f4.6}),
$$
 d\bigl(\TT^{-km}(p,z_0), \TT^{-km}(q_1,z_1)\bigr)
=d\bigl((p, (h^{-k})^4(p,z_0)), (q^{(k)}, (h^{-k})^4(p,z_1)\bigr)\to 0
$$
as $k\to \infty$, and the convergence is exponentially fast.
So by taking the $z$ component, we have that
$|(h^{-k})^4(p,z_0) - (h^{-k+\ell})^4(p,z_1)|\to 0$
converges exponentially fast as $k\to \infty$.
On the other hand, if $z'\not= z''$, then
$|(h^{-k})^4(p,z') - (h^{-k})^4(p,z'')|\to 0$
with a subexponential rate because both
$(h^{-k})^4(p,z')$ and $(h^{-k})^4(p,z'')$ converge to $(p,0)$,
$Dh=\id$ at $(p,0)$, and $h$ is a $C^\infty$ diffeomorphism.
So we get $z_1=(h^{-\ell})^4(p,z_0)$.  This proves (1).

By using the same arguments, and the fact that
$h^4(0, y, 0, z) \ge z$ for any $z\in (0,1)$, we can show
that the $z$-coordinate of $h$ is non-increasing from
$(p,z_5)$ to $(q_2,z_4)$ along $V^s((p,z_5), \TT)$.
That is, $z_4\ge z_5$.  This is (5).

Since the periodic orbit ${\cal C}'$ and the sets
$\cup ^{\infty}_{i=0}G^{-n_1+i}V^s(G^{n_1}(p_1 ), G)\times (0,1)$ and
$\cup^{\infty}_{i=0}G^{n_2-i}V^u(G^{-n_2}(p_2 ), G)\!\times (0,1)$
are unperturbed, we know that the $z$-coordinates
are constant along the stable leaves
$\TT^{-n_1}V^s((G^{n_1}(p_1 ), z), \TT)$
and the unstable leaves $\TT^{n_2}V^u((G^{-n_2}(p_1 ), z), \TT)$.
So we get $z_1=z_2$ and $z_3=z_4$, which are (2) and (4).

Now we prove (3).

Denote by ${\cal C}'(p_1p_2)$ the part of closed orbit ${\cal C}'$
from $p_1$ to $p_2$.
By Lemma~\ref{L7.1}, for any small $\e>0$, we can choose a closed
orbit ${\cal C}'_\e\subset \O_0$ such that a part of
${\cal C}'_\e$ is in the $\e$-neighborhood of ${\cal C}'(p_1p_2)$
in the sense that for any $p'\in {\cal C}'(p_1p_2)$, there is
$p''\in {\cal C}'_\e$ with $d(p',p'')\le \e$.

Consider the map $H=H_{{\cal C}'{\cal C}'_\e}$ from
${\cal C}'(p_1p_2)$ to ${\cal C}'$ given in Section 7.
Starting from $p_1$, we can apply the map consequently
to get a sequence of points $H^{(1)}(p_1), H^{(2)}(p_1), \cdots$, 
where $H^{(i+1)}(p_1)=H(H^{(i)}(p_1))$.
Let $k'$ be such a number that
$H^{(k')}(p_1 )\in {\cal C}'(p_1p_2)$
and $H^{(k'+1)}(p_1 )\notin {\cal C}'(p_1p_2)$.
Clearly, $p_1$ and $p^{(k')}$ can be joint by $4k'$ local stable
and unstable manifolds of $G$ at points on either ${\cal C}'$
or ${\cal C}'_\e$.  Recall that the neighborhoods of both
${\cal C}'\times I$ and ${\cal C}'_\e\times I$ are unperturbed,
and so are the $4k'$ local stable and unstable manifolds.
We know then that $(p_1, z_2)$ is accessible to $(H^{(k')}(p_1), z_2)$
with respect to $\TT$.

Note that $H^{(k')}(p_1 )\in {\cal C}'(p_1p_2)$
and $H^{(k'+1)}(p_1 )\notin {\cal C}'(p_1p_2)$.
By continuity we can find a closer orbit ${\cal O}$
such that $H_{{\cal C}'{\cal O}}(H^{(k')}(p_1 ))=p_2$.
Hence there is $z_3\in (0,1)$ such that $(H^{(k')}(p_1), z_2)$
and $(p_2, z_3)$ are accessible.
We can make $H^{(k')}(p_1)$ arbitrarily close to $p_2$
by taking $\e$ sufficiently small.
Also note that the strong stable and unstable manifolds
change continuously with respect to the diffeomorphism
in the space of partially hyperbolic systems.
If $H^{(k')}(p_1)$ is sufficiently close to $p_2$,
then the $4$ strong stable and unstable
leaves to construct $H_{{\cal C}'{\cal O}}$
are sufficiently short.  Hence, $z_3$ can be made
arbitrarily close to $z_2$.
\qed \enddemo

\section{Hyperbolicity of the Map $Q$}
\setcounter{equation}{0}
\begin{Proposition}\label{P5.1}
There exists $\alpha_0 >0$ such that for any $\alpha\in (0, \alpha_0)$,
$$
\int_N \log \kk_{\alpha 0}(w)dw < \l ^u(S).
$$
\end{Proposition}

\proof
Since $Dh^{(2)}(E^{un}(S))=E^{un}(S)$ for any $w\in N$,
given $\alpha >0$, there exists a unique number $\aa_{\alpha}(w)$
such that the vector
$v_{\alpha}(w)=(1,0,0,\aa_{\alpha}(w))\in E^u (Q_{\alpha 0})$.
Hence $DQ_{\alpha 0}(w)v_{\alpha}(w)=
(\kk_{\alpha}(w), 0,0, \kk_{\alpha}\aa_{\alpha}(w))$
 for some $\kk_{\alpha}> 1$.
The expanding rate of $DQ_{\alpha 0}$ along its unstable direction is
$$
\kk_{\alpha 0}(w)=
\kk_{\alpha}(w)\frac{\sqrt{1+\aa_{\alpha}(Q_{\alpha 0}(w))^2}}
    {\sqrt{1+\aa_{\alpha}(w)^2}}
$$

Let $L_{\alpha}=\int_{N}\log \kk_{\alpha 0}(w)dw=\int_{N}\log
\kk_{\alpha}(w)dw$. The second equality is true because the
map $Q_{\alpha 0}$ preserves the Riemannian volume and therefore
$\int_N \log \bigl(1+\aa_{\alpha}(Q_{\alpha 0}w)^2\bigr) dw
=\int_N \log \bigl(1+\aa_{\alpha}(w)^2\bigr) dw$.
We have that
\begin{eqnarray*}
&Dh^{(2)}_{\alpha}|_{E^{un}}=\left(\begin{array}{cc} A (w,\a )& B(w, \a )
\\C(w,\a ) & D(w, \a ) \end{array}\right)  \\
=&\left(\begin{array}{cc}
r_x \cos \sigma -r\sigma_x \sin\sigma
    & r_t \cos \sigma -r\sigma_t\sin\sigma\\
r_x\sin\sigma +r\sigma_x\cos\sigma
    & r_t\sin\sigma +r\sigma_t\cos\sigma
\end{array}\right)
\end{eqnarray*}
where $r_x =\cos\theta$, $r_t =\sin\theta$ and
$\sigma =\theta+\a\rho (\sqrt{y^2 +z^2 })\rho (r) $
and then $\sigma_x= \frac{\sin\theta}{r}+
\a\rho (\sqrt{y^2 +z^2})\rho '_r (r)\cos\theta$ and
$\sigma_t =\frac{\cos\theta}{r}+
\a\rho (\sqrt{y^2 +z^2})\rho '_r (r)\sin\theta$.

By the same arguments as in the proof of Proposition B.6 in \cite{DHP},
we can get that
\begin{eqnarray*}
\frac{dL_\tau}{d\tau}\Big|_{\tau=0}=0, \qquad
\frac{d^2L_\tau}{d\tau^2}\Big|_{\tau=0}<0.
\end{eqnarray*}
So we can choose $\a_0>0$ so small that $L_\a<\log\l^u(S)$
for any $\a\in (0, \a_0)$.
\qed 
\enddemo


\section{Hyperbolicity of the Map $P$}
\setcounter{equation}{0}

\begin{Proposition}\label{P6.1}
For almost every $w\in N$, $\chi(x, v, P)>0$ $\forall v\in E^{ucn}(x, S)$.
\end{Proposition}

\proof
Denote $\D_0^*=\D_0'\cap \L$.
Then set
$$
U_1=P^{-K}\D_0^*,  \quad
U_2=\D_0\backslash \D_0^*,  \quad
U_3=\D_{k_0}\backslash P^{k_0}\D_0^*,  \quad
U_4=P^{-K}(\D_0\backslash \D_0^*).
$$
Clearly, $\{U_i, i=1,2,3,4\}$ are pairwise disjoint.
Let $U=U_1\cup U_2\cup U_3\cup U_4$.  Note that
$\D^*_0 \supset \C \cap \D_0 '$.
Let $\bar P=P^{\ttau}: U\to U$ be the first return map of $P$,
where $\ttau=\ttau(w)$ is the first return time of $w\in U$.

In the proof below, for any $w\in U$, we always assume
$v\in E^{ucn}(w)$.

If $w\in U_1$, then $\ttau(w)\ge 2K+k_0$.
By Lemma~\ref{L6.2} and (\ref{f3.7}),
\begin{eqnarray}\label{f6.1}
\log\|D{\bar P}_w(v)\|\ge 0.9K\l -0.5\log 2 +\log\|v\|
\ge 0.5K\l  +\log\|v\|.
\end{eqnarray}
By (\ref{f3.12}),
\begin{eqnarray}\label{f6.2}
& \mu U_1 =\mu (\D_0^* )\ge \mu (\Gamma \cap\D_0 ')
\ge \mu (\D_0')-\mu (\Gamma\D\D_0 ')       \nonumber \\
\ge\!\!\! &
\frac{3}{4}\mu\D_0-\mu(\C_0\triangle \D_0)
\ge (0.75-0.05)\mu\D_0=0.7\mu\D_0.
\end{eqnarray}

Note that $\|DP-\id\|\le \d$.  So if $w\in U_2$, then
$\bar P=P^{k_0}$ and
\begin{eqnarray}\label{f6.3}
\log\|D{\bar P}_w(v)\|\ge k_0\log(1-\d)+\log\|v\|.
\end{eqnarray}

Also,
\begin{eqnarray}\label{f6.4}
\mu U_2\!\!\! \!\!\! \!\!\! 
&&=\mu (\D_0 \backslash\D_0^* )
\le\mu (\D_0 \backslash\D_0 ')+\mu (\D_0 '\backslash\D_0^* ) \nonumber  
\le\frac{1}{4}\mu (\D_0 )+\mu (\D_0 ' \backslash(\Gamma\cap\D_0 ')) \\
&& 
\le {1\over 4}\mu\D_0+\mu(\C_0\triangle \D_0)
\le (0.25+0.05)\mu\D_0=0.3\mu\D_0.
\end{eqnarray}

Consider the case that $w\in U_3$.  
The construction of $\{\C_i\}$ and $\{\D_i\}$ implies that
$\ttau(w)\ge K$.
Also, observe that if $P^iw\in \O_3$ for some $i>0$, then 
$P^iw\in \D_j$ for some $0\le j\le k_0-1$.
That is, $P^{i-j}w\in \D_0$.  Hence we have either 
$P^{i-j}w\in U_2\subset U$ or $P^{i-j-K}w\in U_1\subset U$, 
which implies that $\ttau(w)\ge i$, and therefore
the piece of orbit $\{w, \cdots, P^{\ttau(w)-1}w\}$ 
does not intersect $\O_3$.
So we have $\bar P_w (v)=Q^{\ttau(w)}_w (v)$
because $P=Q$ on $N\backslash \O_3$.

Let $\ttau'(w)$ be the smallest positive integer such that
$P^{\ttau'(w)}w\in \L$ for some $0\le \ttau'(w)\le \ttau(w)$,
and let $\ttau'(w)= \tau(w)$ if there is no such integer.
Denote
$$
U_3'=U_3\cap\{w:  \ttau(w)-\ttau'(w)\ge 0.5K\}, \ 
U_3''=U_3\cap\{w:  \ttau(w)-\ttau'(w) < 0.5K\},
$$
and
$$
\hat U_3'=\{P^iw:  w\in U_3', 0\le i < \ttau'(w)\}, \ 
\hat U_3''=\{P^iw:  w\in U_3'', 0\le i < \ttau'(w)\}.
$$

Note that if $n\ge 0.5K$ and $w\in \L$, then
$\|DQ^n_w(v)\|\ge \|v\|$ for any $v\in E^{ucn}(w)$.
Also note that $P=Q$ on $N\backslash \O_3$.

If $w\in U_3'$, then
$$
 \|D{\bar P}_w(v)\|=\|DP^{\ttau(w)}_w(v)\|
=\|DQ^{\ttau(w)-\ttau'(w)}_{P^{\ttau'(w)}(w)}(DP^{\ttau'(w)}_w(v))\|
\ge \|DP^{\ttau'(w)}_w(v))\|.
$$
Hence,
\begin{eqnarray}\label{f6.5}
\qquad\quad  \log\|D{\bar P}_w(v)\|
\ge \log\|P^{\ttau'(w)}_w(v)\|
\ge \sum_{i=0}^{\ttau'(w)-1}\chi_{\hat U_3'}(P^iw)\log(1-\d)+\log\|v\|,
\end{eqnarray}
where $\chi_{\O}(\cdot )$ is the characteristic function of the set $\O$.

If $w\in U_3''$, we denote
$$
\tilde U_3''=\{P^iw: \ w\in U_3'', 0\le i < \ttau(w)\}.
$$
Clearly, we have
\begin{eqnarray}\label{f6.6}
\log \|D{\bar P}_w(v)\|
\ge \sum_{i=0}^{\ttau (w)-1}\chi_{\tilde
  U_3''}(P^iw)\log(1-\d)+\log\|v\|.
\end{eqnarray}

Lastly,  we consider the case that $w\in U_4$.  We have $\ttau(w)=K$
and $\bar P(w)\in U_2$.
We define $\ttau''$ and $U_4'$, $U_4''$, $\hat U_4'$, $\hat U_4''$
and $\tilde U_4''$ in the same way as in the previous case,
and get similar inequalities.  That is,
if $w\in U_4'$, then
\begin{eqnarray}\label{f6.7}
\log\|D{\bar P}_w(v)\|
\ge \sum_{i=0}^{\ttau(w)-1}\chi_{\hat U_4'}(P^iw)\log(1-\d)+\log\|v\|,
\end{eqnarray}
and if $w\in U_4''$, then
\begin{eqnarray}\label{f6.8}
\log\|D{\bar P}_w(v)\|
\ge \sum_{i=0}^{\ttau(w)-1}\chi_{\tilde U_4''}(P^iw)\log(1-\d)+\log\|v\|,
\end{eqnarray}

Clearly, $\hat U_3', \hat U_3'', \hat U_4', \hat U_4''\in \L^c$.
So
\begin{eqnarray}\label{f6.9}
\mu(\hat U_3' \cup \hat U_3''\cup \hat U_4'\cup \hat U_4'')
\le \mu\L^c.
\end{eqnarray}
By the fact that $\ttau'' ,\ttau'\ge 0.5K$ on $\tilde U_3''$
and $\tilde U_4''$ respectively, we have
\begin{eqnarray}\label{f6.10}
\mu(\tilde U_3'' \cup \tilde U_4'')
\le 2\mu(\hat U_3''\cup \hat U_4'')
\le 2\mu\L^c
\end{eqnarray}

Now we estimate the Lyapunov exponent of $v$ at a typical point $w$.
We may assume $w\in U$.
Let $n_i=\ttau_i(w)$ be the $i$th return time of $w$.  Then we have
$$
 \l(v,w,P)
=\lim_{n\to \infty}{1\over n}\log{\|DP^n_w(v)\|\over \|v\|}
=\lim_{i\to \infty}{1\over n_i}
  \log{\|D P^{n_i}_w(v)\|\over \|v\|}
$$
Using the fact that the frequency of the orbit visiting a set $U$
is equal to $\mu U$, and noticing (\ref{f6.1})-(\ref{f6.10}), we have
\begin{eqnarray*}
 &&\lim_{j\to \infty}{1\over n_j}\log\frac{\|DP^{n_j}_wv\|}{\| v\|}
= \lim_{j\to \infty}{1\over n_j}\sum_{i=0}^{j-1}
     \log{\|D\bar P_{\bar P^iw}(D\bar P^i_w(v)\|\over
                              \|D\bar P^i_w(v)\| } \\
&\ge &\!\!\mu U_1\cdot 0.5K\l +\mu U_2\cdot k_0\log(1-\d)
+(\mu\hat U_3' + \mu\hat U_4'+\mu\tilde U_3''+ \mu\tilde U_4')
       \log(1-\d)                \\
&\ge &\!0.7 \mu \D_0\cdot  0.5K\l + 0.3 \mu \D_0 \cdot k_0 \log(1-\d)
+3 \mu \L^c \cdot \log(1-\d)
\end{eqnarray*}

Using (\ref{f3.8}) and (\ref{f3.9}) we conclude
that the right side of the above inequality  is greater than
$$
0.33 \mu\D_0 K\l - 0.003\ttau^2 \l.
$$
Since $\mu\D_0\ge \mu\C - \mu (\C_0 \D\D_0)
\ge 0.08\ttau^2 K^{-1}-0.05 ( 0.08\ttau^2 K^{-1})$
by Lemma~\ref{L6.3} and (\ref{f3.12}),
we conclude finally that $\l (v,w,P)>0$.
\qed 
\enddemo


\begin{Lemma}\label{L6.2}
Let $w\in P^{-K}(\D_0'\cap \L )$.
Then for any $v\in E^{ucn}(w, S)$,
$$
\|D{\bar P}_w(v)\|\ge {\sqrt{2}\over 2}\|v\|e^{0.9K\l}.
$$
\end{Lemma}
\proof
Note that on $\C_{-K,-1}$, $h^{(3)}=\id$, hence $P^K(w)=Q^K(w)$.
Also, since both $Dh^{(1)}$ and $Dh^{(2)}$ preserve the subbundle
$E^{un}(S)$, we have $E^{un}(w, Q)=E^{un}(w, S)$.

Write $v=v^{un}+v^c$, where $v^{un}\in E^{un}(w, Q)$
and $v^c\in E^c(w, Q)$.

We assume first $\disp \|v^c\|\le {\sqrt{2}\over 2}\|v\|$.
Hence $\disp \|v^{un}\|\ge {\sqrt{2}\over 2}\|v\|$.
Since $DQ^K(v^{un})\in E^{un}(Q^Kw, Q)$
and $Q^Kw\in  \Lambda$, we have
$$
\|v^{un}\|=\|DQ^{-K}(DQ^Kv^{un})\|
\le \|DQ^Kv^{un}\|e^{-0.9K\l}.
$$
Hence,
$$
\|DP^Kv\|=\|DQ^Kv\|\ge \|DQ^Kv^{un}\|\ge \|v^{un}\|e^{0.9K\l}
\ge {\sqrt{2}\over 2}\|v\|e^{0.9K\l(Q)}.
$$
Note that at the points $P^Kw, \cdots, P^{K+k_0-1}w$,
the map $Dh^{(3)}$ is a rotation, and
$DQ|_{E^{ucn}(P^iw)}=DS|_{E^{ucn}(P^iw)}$ is noncontracting
for $i=\!K, \cdots, K+k_0-1$.
So $DP^{k_0}|_{E^{ucn}(P^Kw)}$ is noncontracting.
Further, since $\{P^iw\}_{i=K+k_0}^{\ttau}\cap \O_3 =\emptyset$
and $P^{K+k_0}w\in \L'$,
we have that $DP^{\ttau -(K+k_0 )}|_{E^{un}(P^{K+k_0}w)}
=DQ^{\ttau -(K+k_0 )}|_{E^{un}(P^{K+k_0}w)}$ is expanding, and
$DP^{\ttau -(K+k_0 )}|_{E^{ucn}(P^{K+k_0}w)}$ is noncontracting.
So we have
\begin{eqnarray*}
\|D \bar P(v)\|
=&\|DP^{\ttau -(K+k_0 )}_{P^{K+k_0}w} (DP^{K+k_0}_w (v)\|
= \|DQ^{\ttau -(K+k_0 )}_{P^{K+k_0}w} (DP^{K+k_0}_w (v))\|  \\
\ge&\|DP^{K+k_0}_w (v)\|
 \ge\|DP^K _w (v)\| =\|DQ^K _w (v)\|
\ge {\sqrt{2}\over 2}\|v\|e^{0.9K\l}.
\end{eqnarray*}

Now we consider the case that $\disp \|v^c\|\ge {\sqrt{2}\over 2}\|v\|$.

Note that $DQ^K(v^c)\in E^c(Q^{K}w, Q)$.  By the construction of
$h^{(3)}$, we see that $DP^{k_0}_{Q^{K}w}$ rotate the vector in
$E^{cn}(Q^{K}w, S)$ by $\pi/2$.  It means that
$DP^{K+k_0}(v^c)=DP^{k_0}(DQ^K(v^c))\in E^{un}(P^{K+k_0}w, Q)$.
Using the fact that $P^{K+k_0}w\in \L$, we have
$$
\|D \bar P(v^c )\|=\|DP^{\ttau -(K+k_0 )}_{P^{K+k_0}w} (P^{K+k_0}_w (v^c )\|
\ge \|DP^K(DP^{K+k_0}(v^c))\|
$$
$$
 \ge \|DP^{K+k_0}(v^c)\|e^{0.9K\l}
 \ge \|v^c\|e^{0.9K\l}
\ge {\sqrt{2}\over 2}\|v\|e^{0.9K\l}.
$$
This is the result.
\qed \enddemo

\begin{Lemma}\label{L6.3}
$\disp \mu\C 
\ge 0.08\tau^2K^{-1}$.
\end{Lemma}

\proof
Let
$$
\hat \C'=\bigcup_{i=0}^{5\tau^{-1}K-2K-k_0}Q^i\C'.
$$
Since $\tau^{-1}\ge 2$ and $K\ge 2k_0$, we have
$$
{\mu \hat \C'\over \mu \overline \C'}
=\frac{5\tau^{-1}K-2K-k_0 +1}{5\tau^{-1}K-2K-1+K+1}
=1-\frac{K-k_0 +1}{5\tau^{-1}K-K}
\ge 1-0.5\tau.
$$
By (\ref{f3.11}), $\mu \hat \C'=(1-0.5\tau)\mu \overline \C'
\ge 1-0.6\tau$.  Then by (\ref{f3.10}),
$$
\mu (\hat \C\backslash\O)\ge 0.4\tau.
$$

For $w\in \C'$, we denote
$O(w)=\{Q^iw:  i=0, \cdots, (5\tau^{-1}-2)K-k_0\}$,
the piece of orbit that start at $w$ from time $0$ to
$(5\tau^{-1}-2)K-k_0$.
Let
$$
\C_a'=\{O(w): w\in \C', O(w)\cap \O^c\not=\emptyset\}, \
\C_b'=\{O(w): w\in \C', O(w)\cap \O^c=\emptyset\}.
$$
Clearly, $\{\C_a', \C_b'\}$ forms a partition of $\hat\C'$,
and $\C_b'\subset \O$ and therefore by (\ref{f3.10}),
$$
    \mu\C_a'=\mu\hat\C'-\mu\C_b'
\ge \mu\hat\C'-\mu\O
\ge (1-0.6\tau)-(1-\tau)=0.4\tau.
$$
Note that $\C$ consists of exactly one point from each orbit $O(w)$ in
$\C_a$.  We get
$$
\mu \C \ge {\mu \C_a\over (5\tau^{-1}-2)K-k_0+1}
\ge {0.4\tau \over (5\tau^{-1}-2)K}
\ge {0.4\tau \over 5\tau^{-1}K}
=0.08\tau^2K^{-1}.
$$
This is the result.
\qed 
\enddemo

\section{Properties of geodesic flows}
\setcounter{equation}{0}

Let $\PPhi^t: M_0\to M_0$ be the geodesic flow on a compact surface
of a negative constant curvature.
We list some properties of the flow here.

1) $d(\PPhi^tx, x)\le |t|$ for any $t\in \Bbb R$ and $x\in M_0$;

2) $\PPhi^t$ is a uniformly hyperbolic flow, that is, 
there is a decomposition of the tangent bundle into
$$TM_0=E^u\oplus E^s \oplus E^c$$
and a constant $\tilde\kk >1$ such that for any $z\in M_0$,
\begin{eqnarray*}
|D\PPhi^t_z(v)|\ge \tilde\kk^t|v| \qquad v\in E^u_z,  \\
|D\PPhi^{-t}_z(v)|\ge \tilde\kk^{-t}|v| \qquad v\in E^s_z,
\end{eqnarray*}
and $E^c$ is the one dimensional bundle tangent to the flow.

3) The closed orbits are dense in $M_0$.  Moreover,
for any closed orbit ${\cal C}$, both
$W^u({\cal C})=\cup_{z\in {\cal C}}W^u(z)$ and
$W^s({\cal C})=\cup_{z\in {\cal C}}W^s(z)$ are dense in $M_0$.

4) $\PPhi^t$ preserves the Riemannian volume and for any $t\not= 0$,
$\PPhi^t$ is ergodic with respect to the volume,

5) $\PPhi^t$ has the accessibility property.  That is, any two points
can be joint by a piecewise differentiable piecewise nonsingular
path which consists of segments tangent to either $E^u$ or $E^s$.

6) $\PPhi^t$ is topologically conjugate to a symbolic flow that is
a suspension of a subshift finite type with a continuous roof function.
More precisely, there is a symbolic space $\S_A$, two sided left shift
$\s_A: \S_A\to \S_A$, a continuous function $\iota: \S_A\to \Bbb R^+$
and a finite to one map $\pi: \S_A^t\to M_0$ such that
$$
\pi\circ g^t=\bar\s_A^t\circ \pi,
$$
where
$$
\S_A^t=\{(\underline w,t)\in \S_A\times \Bbb R\}/
\{ (\underline w, \iota(\underline w ))=(\s_A(\underline w), 0)\ \},
$$
$$
\bar\s_A^t(\underline w, s)=(\underline w, s+t).
$$
Moreover, $p\in M_0$ is periodic under $\PPhi^t$ if and only if
$\pi^{-1}(p)$ is periodic under $\bar \s_A^t$.
Also, $\S_A$ can be chosen in such a way that the size
of $\pi R_\a$ can be arbitrarily small, where
$ R_\a=\{ w  =\cdots w_{-1} w_0 w_1\cdots: \ w_0 =a\}$
is a cylinder in $\S_A$.  (See \cite{B} and \cite{BR} for more details.)

\begin{Definition}\label{D7.1}
Let ${\cal C}_\e$ and ${\cal C}'$ be two orbits,
$p\in {\cal C}'$, $x\in {\cal C}_\e$ and $d(p,x)$ is small. 
Define the holonomy
$H_{{\cal C}' ,{\cal C}_\e} : B_{{\cal C}'}(p,\d )\rightarrow {\cal
  C}'$ 
where $B_{{\cal C}'}(p,\d )$ is
a $\d$-neighborhood of point $p$ in ${\cal C}'$. $H_{{\cal C}' , 
{\cal C}_\e}
 (x_0 )$ is constructed
the following way (see Figure 2):
\begin{enumerate}
\item $x_1 = W^s (x_0 )\cap W^{uc}(x)$;
\item $x_2 = W^u (x_1 )\cap W^{cs}(x)\subset W^u (x_1 )\cap {\cal C}_\e$;
\item $x_3 = W^s (x_2 )\cap W^{cu}(x_0 )$;
\item $x_4 =W^u (x_3 )\cap W^{cs}(x_0 )\subset W^u (x_3 )\cap {\cal C}'$;
\item $H_{{\cal C}' , {\cal C}_\e} (x_0 )=x_4$.
\end{enumerate}

\begin{figure}
\begin{center}
\psfrag{C'}{$ C'$}
\psfrag{Ce}{$C_\epsilon$}

\psfrag{p}{$p$}

\psfrag{x0}{$x_0$}

\psfrag{x1}{$x_1$}
\psfrag{x2}{$x_2$}

\psfrag{x3}{$x_3$}

\psfrag{x4}{$x_4$}

\psfrag{x}{$x$}

\psfrag{Ws(x0)}{$W^s (x_0 )$}

\psfrag{Wu(x3)}{$W^u (x_3 )$}

\psfrag{Ws(x2)}{$W^s (x_2 )$}

\psfrag{Wu(x1)}{$W^u (x_1)$}

\includegraphics{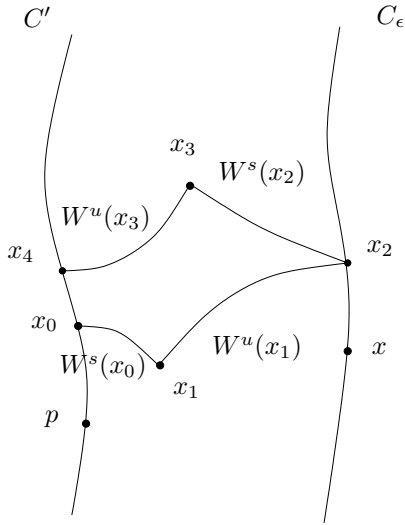}
\end{center}
\caption{Holonomy on a closed orbit.}
\end{figure}
\end{Definition}

This holonomy map was used in \cite{NT}.
Note that in our case for any point $x_0\in {\cal C}'$,
$H_{{\cal C}' ,{\cal C}_\e}(x_0)\neq x_0$
because $g$ is a geodesic flow on a compact surface
of a negative curvature.  

\begin{Lemma}\label{L7.1}
Let ${\cal C}'$ be a given closed  orbits of $\PPhi^t$.
For any $\tau_1\in (0,1)$, there is an open subset $\O_0\subset M_0$ 
containing ${\cal C}'$ with $\mu \O_0\le \tau_1$ and 
there exists  a constant $\e_0>0$ 
such that for any $\e\in (0,\e_0)$,
there is a closed orbit ${\cal C}_\e$ with
$d({\cal C}' ,{\cal C}_\e )<\e$
and $B({\cal C}_\e, \e)\subset \O_0$, where $B({\cal C}_\e, \e)$ 
denotes the $\e$ neighborhood of ${\cal C}_\e$ in $M_0$.
\end{Lemma}

\proof
Let $\s_A: \S_A\to \S_A$ be a symbolic flow that is conjugate to
$g^t: M_0\to M_0$ with the conjugacy $\pi: \S_A\to M_0$.

For a word $W=s_i\cdots s_j$, $i\le j$,
we denote by ${\cal C}(W)$ the periodic orbit of $\s_A$,
and by ${\cal R}(W)$ the cylinder determined by $W$.
Abusing notations we denote by $\pi{\cal C}(W)$ the
corresponding closed orbit of $\PPhi^t$ and by $\pi{\cal R}(W)$
the corresponding set in $M_0$.  More precisely, the latter
means the set
$$
\{(\pi\underline w,t)\in M_0: \underline w\in {\cal R}(W),
0\le t\in \iota(\underline w )\}.
$$

Let $W'$ be the word such that $\pi{\cal C}(W')={\cal C}'$ and
$W_a$ be a subword.
Take another word $W^*$ that generates a periodic orbit
of $\s_A$ and  contains the same subword $W_a$.
Such word $W^*$ exists since the periodic orbits of $\s_A$ are
dense.  We may assume that $W'=W_aW_b$ and $W^*=W_aW_c$
for some word $W_b$ and $W_c$.  It implies that one of the words $W'$
and $W^*$ can be followed by another.

Note that the maximal volume of all sets of the form 
$\pi{\cal R}(W(n))$ converges to  $0$ exponentially fast, 
where $W(n)$ is a $(2n+1)$-word of the form $s_{-n}\cdots s_n$.
For any $\tau _1 >0$, we can take $n>0$ such that
$\mu \bigl(\pi{\cal R}(W(n))\bigr)\le \tau _1 /(4n+2)$ for any 
$(2n+1)$-word $W(n)$.
We may assume $n>|W'|, |W^*|$, where $|W|$ denotes the length
of $W$.   Note that for any integer $k\ge n/|W'|$, there are
at most $4n+2$ different $(2n+1)$-words of the form $W(n)$
in the orbit of $\underline w={\cal C}((W')^k W^*)$,
where $W^k$ is the word consists $k$ consecutive $W$.
Let $W_1, \cdots, W_j$, $j\le 4n+2$, denote these words.

Set 
$$
\hat \O_0=\bigcup_{i=1}^j \pi{\cal R}(W_i).
$$
Clearly, $\mu\hat \O_0 \le \tau_1 /2$.  Choose $\e_0>0$
such that $\mu B(\hat \O_0, \e_0)\le \tau_1$.
Then we set $\O_0=B(\hat \O_0, \e_0)$.

Let $\e\in (0,\e_0)$.  We take a word of the form 
$(W')^k W^*$ for some large $k$, and then 
take ${\cal C}_\e=\pi({\cal C}((W')^k W^*)$.
Clearly, ${\cal C}_\e\subset \hat \O_0$, 
and therefore $B({\cal C}_\e, \e)\subset \O_0$ 
by the choice of $\O_0$.   
Also if $k$ is large enough, then the distance 
between ${\cal C}_\e$ and ${\cal C}'$ can be made
arbitrarily small.  Hence we have  
${\cal C}'\subset B({\cal C}_\e, \e)$. 
\qed

\begin {Lemma}\label{L6.4}
For any $\c>0$, $\tau_2\in (0,1)$ and $C>0$, 
we can find a constant $\t_0>0$ such that for any number 
$\t\in [0, \t_0]$, any set of the form
$\D=\D_{ss's''}=B^u(x, s')\times B^s(y, s'')\times B^{cn}((t,z), s)$,
where $s', s''\ge s$ and $s', s''\le C$, there exists a set 
$\D^{(0)}$ and a map $h: N\to N$ with the following properties:
\begin{itemize}
\item[{\rm (a)}]  $h=T_\t$ on $\D^{(0)}$, and $h=\id$ on $\D^c$;
\item[{\rm (b)}]  $\mu \D^{(0)}/\mu\D\ge \tau_2$;
\item[{\rm (c)}]  $\|h-\id\|\le \c$,
\end{itemize}
where $T_\t$ is a rotation given by
$$
T_\t(x,y,t,z)=(x, y, t\cos\t-z\sin\t, t\sin\t+z\cos\t).
$$
\end{Lemma}

\proof
Take $\kappa>0$ such that
$\mu \D_{1-\kappa, 1-\kappa, 1-\kappa}/\mu \D_{111}\ge \tau_2$.
Hence, for any $r>0$, $r', r''>r$, we have
$\mu \D_{r-\kappa r, r'-\kappa r, r''-\kappa r}/\mu \D_{rr'r''}\ge \tau_2$,
since $r'/(r'-\kappa r)$ and $r''/(r''-\kappa r)$ are increasing.

Take a family of $C^\infty$ functions
$\zeta_r=\zeta_r(s): \Bbb R^+\to \Bbb R^+$, for $r\ge 1$
such that
\begin{enumerate}
\item[1] $\zeta_1(s)=1$ if $s\in [0, 1-\kappa]$ and
$\zeta_1(s)=0$ if $s\ge 1$;
\item[2] $\zeta_r(s)=1$ if $s\in [0, r-1)$ and
$\zeta_r(s)=\zeta_1(s-r+1)$ if $s\ge r-1$.
\end{enumerate}
Clearly, $\zeta_r$ have that same $C^\infty$ norm for all $r\ge 1$.

Take coordinate system $w=(x,y,t,z)$ as in (\ref{f3.0}).
Then we define $h:\D\to \D$ by
$$
h(w)=T_{\t(s,s',s'')}(w), \qquad
\t(s,s',s'')=\t\zeta_{s'/s}(x/s)\zeta_{s''/s}(y/s)
  \zeta_1(\sqrt{t^2+z^2}/s).
$$

By the construction, we see that $h$ satisfies (a) and (b).
For (c),
note that if $\t=0$, then $h=\id$,
and note that the $C^1$ norm of $h$ change smoothly with $\t$,
we get the result. 
\qed \enddemo

\section*{Acknowledgments}

The authors would like to thank Professor Yakov Pesin for suggesting 
this topic and for many important discussions.  They would also
like to thank Professor M. Brin for a useful discussion. 

Part of the work was done when the first author was in the Department 
of Mathematics, University of Southern California, Los Angeles, CA 90089.
This work was partially supported by NSF under grants DMS-0196234 
and DMS-0240097.


\end{document}